\documentclass{amsart}
\usepackage{amssymb}

\newcommand{\arx}[1]{\texttt{http://arxiv.org/abs/#1}}

\newcommand{\Weiss}{
\author{Tomasz Weiss}
\address{Institute of Mathematics, Akademia Podlaska
08-119 Siedlce, Poland}
\email{tomaszweiss@go2.pl}
}

\newcommand{\Tsaban}{
\author{Boaz Tsaban}
\thanks{The second author is partially supported by the Golda Meir Fund and
the Edmund Landau Center for Research in Mathematical Analysis and Related Areas,
sponsored by the Minerva Foundation (Germany).}
\address{Einstein Institute of Mathematics, Hebrew University of Jerusalem,
Givat Ram, Jerusalem 91904, Israel}
\email{tsaban@math.huji.ac.il, http://www.cs.biu.ac.il/\~{}tsaban}
}

\newcommand{\CH}{the Continuum Hypothesis}
\long\def\forget#1\forgotten{}

\renewcommand{\pi}{pseudointersection}
\newcommand{\diam}{\op{diam}}

\newcommand{\cU}{\mathcal{U}}

\newcommand{\fU}{\mathfrak{U}}
\newcommand{\fV}{\mathfrak{V}}
\newcommand{\seq}[1]{\{#1\}_{n\in\N}}

\newcommand{\op}{\operatorname}
\newcommand{\hdim}{\op{dim}}

\newcommand{\Cantor}{{\{0,1\}^\N}}

\newcommand{\threeCantor}{{\{-1,0,1\}^\N}}
\newcommand{\Cal}{\mathcal}

\newcommand{\B}{{\Cal B}}

\newcommand{\BG}{\B_\Gamma}

\newcommand{\BO}{\B_\Omega}

\newcommand{\Tau}{\mathrm{T}}

\newcommand{\cF}{{\Cal F}}

\newcommand{\N}{\mathbb{N}}

\renewcommand{\O}{\Cal O}

\newcommand{\Q}{\mathbb{Q}}
\newcommand{\R}{\mathbb{R}}

\newcommand{\Union}{\bigcup}

\renewcommand{\c}{{\mathfrak c}}

\renewcommand{\i}{\item}

\newcommand{\oo}{\infty}
\newcommand{\p}{{\mathfrak p}}

\newcommand{\x}{\times}

\newcommand{\nin}{\not\in}

\newcommand{\sbst}{\subseteq}

\newcommand{\sm}{\setminus}

\newcommand{\as}{\sbst^*}
\renewcommand{\|}{\upharpoonright}
\renewcommand{\(}{\left(}
\renewcommand{\)}{\right)}

\newcommand{\impl}{\to}

\newtheorem{thm}{Theorem} 
\newtheorem{prop}[thm]{Proposition}
\newtheorem{prob}[thm]{Problem}

\newtheorem{lem}[thm]{Lemma}

\theoremstyle{definition}

\theoremstyle{remark}

\newcommand{\be}{\begin{enumerate}}
\newcommand{\ee}{\end{enumerate}}
\newcommand{\bi}{\begin{itemize}}
\newcommand{\ei}{\end{itemize}}

\newcommand{\sone}{{\sf S}_1}

\newcommand{\lft}[2]{\mathopen\ifcase#1{}\oo\or
                        \big#2\or\Big#2\else\oo\fi}
\newcommand{\rgt}[2]{\mathclose\ifcase#1{}\oo\or
                        \big#2\or\Big#2\else\oo\fi}

\title{Topological diagonalizations and Hausdorff dimension}

\Weiss
\Tsaban



\begin{document}
\begin{abstract}
The Hausdorff dimension of a product
$X\x Y$ can be strictly greater than
that of $Y$, even when the Hausdorff dimension of
$X$ is zero. But when $X$ is countable, the Hausdorff dimensions
of $Y$ and $X\x Y$ are the same.
Diagonalizations of covers define a natural hierarchy of properties which
are weaker than ``being countable'' and stronger than
``having Hausdorff dimension zero''.
Fremlin asked whether it is enough for $X$ to have the strongest
property in this hierarchy (namely, being a $\gamma$-set)
in order to assure that the Hausdorff dimensions
of $Y$ and $X\x Y$ are the same.

We give a negative answer: Assuming \CH{},
there exists a $\gamma$-set $X\sbst\R$
and a set $Y\sbst\R$ with Hausdorff dimension zero, such that the Hausdorff
dimension of $X+Y$ (a Lipschitz image of $X\x Y$) is maximal, that is, $1$.
However, we show that for the notion of a
\emph{strong} $\gamma$-set the answer is positive.
Some related problems remain open.
\end{abstract}
\maketitle

\section{Introduction}
The Hausdorff dimension of a subset of $\R^k$ is a derivative of
the notion of Hausdorff \emph{measures} \cite{Fa}. However, for
our purposes it will be more convenient to use the following
equivalent definition. Denote the diameter of a subset $A$ of
$\R^k$ by $\diam(A)$. The \emph{Hausdorff dimension}
of a set $X\sbst\R^k$, $\hdim(X)$, is the infimum of all positive $\delta$ such
that for each positive $\epsilon$ there exists a cover
$\seq{I_n}$ of $X$ with
$$\sum_{n\in\N}\diam(I_n)^\delta<\epsilon.$$
From the many properties of Hausdorff dimension, we will need the
following easy ones.
\goodbreak
\begin{lem}
~\be
\i If $X\sbst Y\sbst\R^k$, then $\hdim(X)\le\hdim(Y)$.
\i Assume that $X_1,X_2,\dots$ are subsets of $\R^k$ such that
$\hdim(X_n)=\delta$ for each $n$. Then $\hdim(\Union_n X_n)=\delta$.
\i Assume that $X\sbst\R^k$ and $Y\sbst\R^m$ is such that there exists
a Lipschitz surjection $\phi:X\to Y$. Then $\hdim(X)\ge\hdim(Y)$.
\i For each $X\sbst\R^k$ and $Y\sbst\R^m$, $\hdim(X\x Y)\ge \hdim(X)+\hdim(Y)$.
\ee
\end{lem}
Equality need not hold in item (4) of the last lemma.
In particular, one can construct
a set $X$ with Hausdorff dimension zero and a set $Y$ such that
$\hdim(X\x Y)>\hdim(Y)$.
On the other hand, when $X$ is countable, $X\x Y$ is a union of
countably many copies of $Y$, and therefore
\begin{equation}\label{equaldim}
\hdim(X\x Y) = \hdim(Y).
\end{equation}
Having Hausdorff dimension zero can be thought of as a notion of smallness.
Being countable is another notion of smallness, and we know that
the first notion is not enough restrictive in order to have Equation \ref{equaldim}
hold, but the second is.

Notions of smallness for sets of real numbers have a long history and many applications --
see, e.g., \cite{MilSpec}.
We will consider some notions which are
weaker than being countable and stronger than having Hausdorff
dimension zero.

According to Borel \cite{BOREL},
a set $X\sbst\R^k$ has \emph{strong measure zero} if
for each sequence of positive reals $\seq{\epsilon_n}$,
there exists a cover
$\seq{I_n}$ of $X$ such that $\diam(I_n)<\epsilon_n$ for all $n$.
Clearly strong measure zero implies Hausdorff dimension zero.
It does not require any special assumptions in order to see
that the converse is false.
A perfect set can be mapped onto the unit interval by a uniformly
continuous function and therefore cannot have strong measure zero.

\begin{prop}[folklore]\label{smallCantor}
There exists a perfect set of reals $X$ with Hausdorff dimension zero.
\end{prop}
\begin{proof}
For $0<\lambda<1$, denote by $C(\lambda)$ the Cantor set obtained
by starting with the unit interval, and at each step removing
from the middle of each interval a subinterval of size
$\lambda$ times the size of the interval
(So that $C(1/3)$ is the canonical middle-third Cantor set, which has
Hausdorff dimension $\log 2/\log 3$.)
It is easy to see that if $\lambda_n\nearrow 1$, then
$\hdim(C(\lambda_n))\searrow 0$.

Thus, define a special Cantor set $C(\seq{\lambda_n})$
by starting with the unit interval, and at step $n$ removing
from the middle of each interval a subinterval of size
$\lambda_n$ times the size of the interval.
For each $n$, $C(\seq{\lambda_n})$ is contained in
a union of $2^n$ (shrunk) copies of $C(\lambda_n)$,
and therefore $\hdim(C(\seq{\lambda_n}))\le\hdim(C(\lambda_n))$.
\end{proof}

As every countable set has strong measure zero, the latter notion can be thought
of an ``approximation'' of countability. In fact, Borel conjectured
in \cite{BOREL} that every strong measure zero set is countable, and
it turns out that the usual axioms of mathematics (ZFC) are not strong enough
to prove or disprove this conjecture: Assuming \CH{} there exists
an uncountable strong measure zero set (namely, a Luzin set),
but Laver \cite{LAVER} proved that one cannot prove
the existence of such an object
from the usual axioms of mathematics.

The property of strong measure zero (which depends on the metric)
has a natural topological counterpart.
A topological space $X$ has \emph{Rothberger's property $C''$} \cite{ROTH41}
if for each sequence $\seq{\cU_n}$ of covers of $X$
there is a sequence $\seq{U_n}$ such that for each $n$ $U_n\in\cU_n$, and $\seq{U_n}$
is a cover of $X$.
Using Scheepers' notation \cite{coc1},
this property is a particular instance of the following selection hypothesis
(where $\fU$ and $\fV$ are any collections of covers of $X$):
\bi
\item[$\sone(\fU,\fV)$:]
For each sequence $\seq{\cU_n}$ of members of $\fU$, there is a
sequence $\seq{U_n}$ such that $U_n\in\cU_n$ for each $n$, and
$\seq{U_n}\in\fV$.
\ei
Let $\O$ denote the
collection of all open covers of $X$.
Then the property considered by Rothberger is $\sone(\O,\O)$.
Fremlin and Miller \cite{FM} proved that a set $X\sbst\R^k$ satisfies
$\sone(\O,\O)$ if, and only if, $X$ has strong measure zero with respect to
each metric which generates the standard topology on $\R^k$.

But even Rothberger's property for $X$ is not strong enough to have
Equation \ref{equaldim} hold: It is well-known that every Luzin set satisfies
Rothberger's property (and, in particular, has Hausdorff dimension zero).
\begin{lem}
The mapping $(x,y)\mapsto x+y$ from $\R^2$ to $\R$ is Lipschitz.
\end{lem}
\begin{proof}
Observe that for nonnegative
reals $a$ and $b$, $(a-b)^2\ge 0$ and therefore $a^2+b^2\ge 2ab$.
Consequently,
$$a+b = \sqrt{a^2+2ab+b^2}\le\sqrt{2(a^2+b^2)} = \sqrt{2}\sqrt{a^2+b^2}.$$
Thus, $|(x_1+y_1)-(x_2+y_2)|\le\sqrt{2}\sqrt{(x_1-x_2)^2+(y_1-y_2)^2}$
for all $(x_1,y_1),(x_2,y_2)\in\R^2$.
\end{proof}
Assuming \CH{}, there exists a Luzin set
$L\sbst\R$ such that $L+L$, a Lipschitz image of $L\x L$, is equal to $\R$
\cite{coc2}.

We therefore consider some stronger properties.
An open cover $\cU$ of $X$ is an
\emph{$\omega$-cover} of $X$ if each finite subset of
$X$ is contained in some member of the cover, but $X$ is not
contained in any member of $\cU$.
$\cU$ is a $\gamma$-cover of $X$ if it is infinite, and
each element of $X$ belongs to all but finitely many
members of $\cU$.
Let $\Omega$
and $\Gamma$ denote the
collections of open $\omega$-covers
and $\gamma$-covers of $X$,
respectively. Then
$\Gamma\sbst\Omega\sbst\O$,
and these three classes of covers introduce $9$ properties of the form $\sone(\fU,\fV)$.
If we remove the trivial ones and check for equivalences \cite{coc2, strongdiags},
then it turns out that only six of these properties are really distinct,
and only three of them imply Hausdorff dimension zero:
$$\sone(\Omega,\Gamma)\impl\sone(\Omega,\Omega)\impl\sone(\O,\O).$$
The properties $\sone(\Omega,\Gamma)$ and $\sone(\Omega,\Omega)$
were also studied before.
$\sone(\Omega,\Omega)$ was studied by Sakai \cite{Sakai},
and $\sone(\Omega,\Gamma)$ was studied by Gerlits and Nagy
in \cite{GN}:
A topological space $X$ is a \emph{$\gamma$-set}
if each $\omega$-cover of $X$ contains a $\gamma$-cover of $X$.
Gerlits and Nagy proved that $X$ is a $\gamma$-set if, and only
if, $X$ satisfies $\sone(\Omega,\Gamma)$.
It is not difficult to see that every countable space is a $\gamma$-set.
But this property is not trivial:
Assuming \CH{}, there exist uncountable $\gamma$-sets \cite{GM}.

$\sone(\Omega,\Omega)$ is closed under taking finite
powers \cite{coc2}, thus the
Luzin set we used to see that Equation \ref{equaldim} need not hold when
$X$ satisfies $\sone(\O,\O)$ does not rule out that possibility that this Equation
holds when $X$ satisfies $\sone(\Omega,\Omega)$.
However, in \cite{huremen2} it is shown that assuming \CH{}, there exist
Luzin sets $L_0$ and $L_1$ satisfying $\sone(\Omega,\Omega)$, such that
$L_0+L_1=\R$. Thus, the only remaining candidate for a nontrivial
property of $X$ where Equation \ref{equaldim} holds is
$\sone(\Omega,\Gamma)$ ($\gamma$-sets).
Fremlin (personal communication) asked whether Equation \ref{equaldim} is indeed
provable in this case.
We give a negative answer, but
show that for a yet stricter (but nontrivial) property which was
considered in the literature, the answer is positive.

The notion of a strong $\gamma$-set was introduced in \cite{GM}.
However, we will adopt the following simple characterization
from \cite{strongdiags} as our formal definition.
Assume that $\seq{\fU_n}$ is a sequence of collections of covers of a space $X$,
and that $\fV$ is a collection of covers of $X$.
Define the following selection hypothesis.
\begin{itemize}
\item[$\sone(\seq{\fU_n},\fV)$:]
For each sequence $\seq{\cU_n}$ where $\cU_n\in\fU_n$ for each $n$,
there is a sequence
$\seq{U_n}$ such that $U_n\in\cU_n$ for each $n$, and $\seq{U_n}\in\fV$.
\end{itemize}
A cover $\cU$ of a space $X$ is an \emph{$n$-cover} if each
$n$-element subset of $X$ is contained in some member of $\cU$.
For each $n$ denote by $\O_n$ the collection of all
open $n$-covers of a space $X$. Then $X$ is a \emph{strong $\gamma$-set} if
$X$ satisfies $\sone(\seq{\O_n},\Gamma)$.

In most cases $\sone(\seq{\O_n},\fV)$ is equivalent to $\sone(\Omega,\fV)$ \cite{strongdiags},
but not in the case $\fV=\Gamma$: It is known that for a strong $\gamma$-set $G\sbst\Cantor$
and each $A\sbst\Cantor$ of measure zero, $G\oplus A$ has measure zero too \cite{GM}; this can
be contrasted with Theorem \ref{BRthm} below.
In Section \ref{sg} we show that Equation \ref{equaldim} is provable in the case that $X$ is
a strong $\gamma$-set, establishing another difference between the notions of $\gamma$-sets
and strong $\gamma$-sets, and giving a positive answer to Fremlin's question under a stronger
assumption on $X$.

\section{The product of a $\gamma$-set and a set of Hausdorff dimension zero}

\begin{thm}\label{gammaprod}
Assuming \CH{} (or just $\p=\c$),
there exist a $\gamma$-set $X\sbst\R$ and a set $Y\sbst\R$ with
Hausdorff dimension zero such that the Hausdorff dimension of
the algebraic sum
$$X+Y = \{ x + y : x\in X, y\in Y\}$$
(a Lipschitz image of $X\x Y$ in $\R$) is $1$.
In particular, $\dim(X\x Y)\ge 1$.
\end{thm}

Our theorem will follow from the following related theorem.
This theorem involves the \emph{Cantor space} $\Cantor$ of infinite
binary sequences. The Cantor space is equipped with the product topology
and with the product measure.

\begin{thm}[Bartoszy{\'n}ski and Rec{\l}aw \cite{BR}]\label{BRthm}
Assume \CH{} (or just $\p=\c$).
Fix an increasing sequence $\seq{k_n}$ of natural numbers, and for each
$n$ define
$$A_n=\{f\in \Cantor : f\| [k_n,k_{n+1})\equiv 0\}.$$
If the set
$$A = \bigcap_{m\in \N }\Union_{n\ge m } A_n$$
has measure zero, then
there exists a $\gamma$-set  $G\sbst\Cantor$
such that the algebraic sum $G\oplus A$ is equal to $\Cantor$ (where
where $\oplus$ denotes the modulo $2$ coordinatewise addition).
\end{thm}
Observe that the assumption in Theorem \ref{BRthm} holds whenever
$\sum_n 2^{-(k_{n+1}-k_n)}$ converges.

\begin{lem}\label{Hdim0}
There exists an increasing sequence of natural
numbers $\seq{k_n}$ such that
$\sum_n 2^{-(k_{n+1}-k_n)}$ converges,
and such that for the sequence
$\seq{B_n}$ defined by
$$B_n = \left\{\sum_{i\in\N}\frac{f(i)}{2^{i+1}} : f\in\threeCantor\mbox{ and }
f\|[k_n,k_{n+1})\equiv 0\right\}$$
for each $n$,
the set
$$Y=\bigcap_{m\in \omega }\bigcup_{n\ge m} B_n$$
has Hausdorff dimension zero.
\end{lem}
\begin{proof}
Fix a sequence $p_n$ of positive reals which converges to $0$.
Let $k_0=0$. Given $k_n$ find $k_{n+1}$ satisfying
$$3^{k_n}\cdot \frac{1}{2^{p_n(k_{n+1}-2)}}\le \frac{1}{2^n}.$$
Clearly, every $B_n$ is contained in a union of $3^{k_n}$ intervals
such that each of the intervals has diameter
$1/2^{k_{n+1}-2}$.
For each positive $\delta$ and $\epsilon$,
choose $m$ such that $\sum_{n\ge m}1/2^n<\epsilon$
and such that $p_n<\delta$ for all $n\ge m$.
Now, $Y$ is a subset of $\Union_{n\ge m}B_n$, and
$$\sum_{n\ge m}3^{k_n}\(\frac{1}{2^{k_{n+1}-2}}\)^\delta <
\sum_{n\ge m}3^{k_n}\(\frac{1}{2^{k_{n+1}-2}}\)^{p_n} <
\sum_{n\ge m} \frac{1}{2^n}<\epsilon.$$
Thus, the Hausdorff dimension of $Y$ is zero.
\end{proof}

The following lemma concludes the proof of Theorem \ref{gammaprod}.

\begin{lem}
There exists a $\gamma$-set $X\sbst\R$ and a set $Y\sbst\R$ with
Hausdorff dimension zero such that $X+Y=\R$. In particular,
$\dim(X+Y)=1$.
\end{lem}
\begin{proof}
Choose a sequence $\seq{k_n}$ and a set $Y$ as in Lemma \ref{Hdim0}.
Then $\sum_n 2^{-(k_{n+1}-k_n)}$ converges, and
the corresponding set $A$ defined
in Theorem \ref{BRthm} has measure zero.
Thus, there exists a $\gamma$-set $G$ such that $G\oplus A = \Cantor$.

Define $\Phi:\Cantor\to\R$ by
$$\Phi(f)=\sum_{i\in\N}\frac{f(i)}{2^{i+1}}.$$
As $\Phi$ is continuous, $X = \Phi[G]$ is a $\gamma$-set of reals.
Assume that $z$ is a member of the interval $[0,1]$, let $f\in\Cantor$
be such that $z=\sum_i f(i)/2^{i+1}$. Then
$f = g \oplus a$ for appropriate $g\in G$ and $a\in A$.
Define $h\in\threeCantor$ by $h(i) = f(i)-g(i)$.
For infinitely many $n$, $a\|[k_n,k_{n+1})\equiv 0$ and therefore
$f\|[k_n,k_{n+1})\equiv g\|[k_n,k_{n+1})$, that is,
$h\|[k_n,k_{n+1})\equiv 0$ for infinitely many $n$. Thus, $y = \sum_i h(i)/2^{i+1}\in Y$,
and for $x=\Phi(g)$,
$$x + y = \sum_{i\in\N}\frac{g(i)}{2^{i+1}} + \sum_{i\in\N}\frac{h(i)}{2^{i+1}} =
\sum_{i\in\N}\frac{g(i)+h(i)}{2^{i+1}} = \sum_{i\in\N}\frac{f(i)}{2^{i+1}} = z.$$
This shows that $[0,1]\sbst X+Y$. Consequently,
$X+(Y+\Q) = (X+Y)+\Q = \R$. Now, observe that $Y+\Q$ has Hausdorff dimension zero since
$Y$ has.
\end{proof}

\section{The product of a strong $\gamma$-set and a set of Hausdorff dimension zero}\label{sg}

\begin{thm}\label{SGprod}
Assume that $X\sbst\R^k$ is a strong $\gamma$-set.
Then for each $Y\sbst\R^l$, $\hdim(X\x Y)=\hdim(Y)$.
\end{thm}
\begin{proof}
The proof for this is similar to that of Theorem 7 in \cite{GM}.

It is enough to show that $\hdim(X\x Y)\le\hdim(Y)$.
\begin{lem}\label{largeHcover}
Assume that $Y\sbst\R^l$ is such that $\hdim(Y)<\delta$.
Then for each positive $\epsilon$ there exists a \emph{large}
cover $\seq{I_n}$ of $Y$ (i.e., such that each $y\in Y$ is a
member of infinitely many sets $I_n$) such that
$\sum_n\diam(I_n)^\delta<\epsilon$.
\end{lem}
\begin{proof}
For each $m$ choose a cover $\seq{I^m_n}$ of $Y$
such that $\sum_n\diam(I^m_n)^\delta<\epsilon/2^m$.
Then $\{I^m_n : m,n\in\N\}$ is a large cover of $Y$,
and $\sum_{m,n}\diam(I^m_n)^\delta<\sum_n\epsilon/2^m=\epsilon$.
\end{proof}
\begin{lem}\label{goodHcover}
Assume that $Y\sbst\R^l$ is such that $\hdim(Y)<\delta$.
Then for each sequence $\seq{\epsilon_n}$ of positive reals
there exists a large cover $\seq{A_n}$ of $Y$
such that for each $n$ $A_n$ is a union of finitely many sets, $I^n_1,\dots,I^n_{m_n}$,
such that $\sum_j \diam(I^n_j)^\delta<\epsilon_n$.
\end{lem}
\begin{proof}
Assume that $\seq{\epsilon_n}$ is a sequence of positive reals.
By Lemma \ref{largeHcover}, there exists
a large cover $\seq{I_n}$ of $Y$ such that
$\sum_n\diam(I_n)^\delta<\epsilon_1$.
For each $n$ let
$k_n = \min\{ m : \sum_{j\ge m}\diam(I_j)^\delta<\epsilon_n\}$.
Take
$$A_n = \Union_{j=k_n}^{k_{n+1}-1}I_j.$$
\end{proof}
Fix $\delta>\hdim(Y)$ and $\epsilon>0$.
Choose a sequence $\seq{\epsilon_n}$ of positive reals such that
$\sum_n 2n\epsilon_n<\epsilon$, and use Lemma \ref{goodHcover}
to get the corresponding large cover $\seq{A_n}$.

For each $n$ we define an $n$-cover $\cU_n$ of $X$ as follows.
Let $F$ be an $n$-element subset of $X$.
For each $x\in F$,
find an open interval $I_x$ such that $x\in I_x$
and
$$\sum_{j=1}^{m_n} \diam(I_x\x I^n_j)^\delta<2\epsilon_n.$$
Let $U_F=\Union_{x\in F}I_x$.
Set
$$\cU_n = \{U_F : F\mbox{ is an $n$-element subset of }X\}.$$
As $X$ is a strong $\gamma$-set, there exist elements
$U_{F_n}\in\cU_n$, $n\in\N$, such that $\seq{U_{F_n}}$ is a $\gamma$-cover of $X$.
Consequently,
$$X\x Y \sbst\Union_{n\in\N} (U_{F_n}\x A_n)\sbst
\Union_{n\in\N}\Union_{x\in F_n}\Union_{j=1}^{m_n} I_x\x I^n_j$$
and
$$\sum_{n\in\N}\sum_{x\in F_n}\sum_{j=1}^{m_n}\diam(I_x\x I^n_j)^\delta<
\sum_n n\cdot 2\epsilon_n<\epsilon.$$
\end{proof}

\section{Open problems}
There are ways to strengthen the notion of $\gamma$-sets other than moving
to strong $\gamma$-sets.
Let $\BO$ and $\BG$ denote the collections of \emph{countable Borel}
$\omega$-covers and $\gamma$-covers of $X$, respectively.
As every open $\omega$-cover of a set of reals contains a countable
$\omega$-subcover \cite{coc2}, we have that $\Omega\sbst\BO$ and therefore
$\sone(\BO,\BG)$ implies $\sone(\Omega,\Gamma)$. The converse is not
true \cite{CBC}.

\begin{prob}
Assume that $X\sbst\R$ satisfies $\sone(\BO,\BG)$.
Is it true that for each $Y\sbst\R$, $\hdim(X\x Y)=\hdim(Y)$?
\end{prob}
We conjecture that assuming \CH{},
the answer to this problem is negative.
We therefore introduce the following problem.
For infinite sets of natural numbers $A,B$, we write $A\as B$
if $A\sm B$ is finite. Assume that $\cF$ is a family of infinite
sets of natural numbers.
A set $P$ is a \emph{\pi{}}
of $\cF$ if it is infinite, and for each $B\in\cF$, $A\as B$.
$\cF$ is \emph{centered} if each finite subcollection of
$\cF$ has a \pi{}.
Let $\p$ denote the minimal cardinality of a centered family
which does not have a \pi{}.
In \cite{CBC} it is proved that
$\p$ is also the minimal cardinality of a set of reals which
does not satisfy $\sone(\BO,\BG)$.

\begin{prob}
Assume that the cardinality of $X$ is smaller than $\p$.
Is it true that for each $Y\sbst\R$, $\hdim(X\x Y)=\hdim(Y)$?
\end{prob}

Another interesting open problem involves the following
notion \cite{tau, tautau}.
A cover $\cU$ of $X$ is a \emph{$\tau$-cover} of $X$
if it is a large cover, and for each $x,y\in X$,
one of the sets $\{U\in\cU : x\in U\mbox{ and }y\nin U\}$
or $\{U\in\cU : y\in U\mbox{ and }x\nin U\}$ is finite.
Let $\Tau$ denote the collection of open $\tau$-covers of $X$.
Then $\Gamma\sbst\Tau\sbst\Omega$, therefore
$\sone(\seq{\O_n},\Gamma)$ implies $\sone(\seq{\O_n},\Tau)$.

\begin{prob}\label{probOnT}
Assume that $X\sbst\R$ satisfies $\sone(\seq{\O_n},\Tau)$.
Is it true that for each $Y\sbst\R$, $\hdim(X\x Y)=\hdim(Y)$?
\end{prob}

It is conjectured that $\sone(\seq{\O_n},\Tau)$
is strictly stronger than $\sone(\Omega,\Tau)$ \cite{strongdiags}.
If this conjecture is false, then the results in this paper imply
a negative answer to Problem \ref{probOnT}.

Another type of problems is the following:
We have seen that the assumption that $X$ is a $\gamma$-set and
$Y$ has Hausdorff dimension zero is not enough in order to
prove that $X\x Y$ has Hausdorff dimension zero.
We also saw that if $X$ satisfies a stronger property
(strong $\gamma$-set), then $\dim(X\x Y)=\dim(Y)$ for all $Y$.
Another approach to get a positive answer would be to strengthen
the assumption on $Y$ rather than $X$.

If we assume that $Y$ has strong measure zero, then
a positive answer follows from a result of Scheepers
\cite{smzpow} (see also \cite{prods}), asserting that
if $X$ is a strong measure zero metric space which
also has the Hurewicz property, then for each strong
measure zero metric space $Y$, $X\times Y$ has strong
measure zero.
Indeed, if $X$ is a $\gamma$-set then it has the
required properties.

Finally, the following question of Krawczyk remains open.

\begin{prob}
Is it consistent (relative to ZFC) that
there are uncountable $\gamma$-sets
but for each $\gamma$-set $X$ and each set $Y$,
$\dim(X\x Y)=\dim(Y)$?
\end{prob}

\end{document}